\documentclass{article}
\usepackage{graphicx}
\usepackage{amsmath,amssymb,amsfonts}

\newcommand{\C}{\mathcal{C}}
\newcommand{\E}{\mathcal{E}}
\newcommand{\U}{\mathcal{U}}
\newcommand{\D}{\mathcal{D}}
\newcommand{\DD}{\mathfrak{D}}
\newcommand{\V}{\mathcal{V}}
\newcommand{\Q}{\mathcal{Q}}
\newcommand{\B}{\mathcal{B}}
\newcommand{\T}{\mathcal{T}}

\newcommand{\ar}{%
  \mathrel{\mathchoice
    {\scalebox{0.5}[1]{$\displaystyle\rightarrow$}}
    {\scalebox{0.5}[1]{$\textstyle\rightarrow$}}
    {\scalebox{0.5}[1]{$\scriptstyle\rightarrow$}}
    {\scalebox{0.5}[1]{$\scriptscriptstyle\rightarrow$}}}
}

\newcommand{\bool}{o}
\newcommand{\Const}{\C\textit{onst}}
\newcommand{\Var}{\V\textit{ar}}
\newcommand{\iiota}{\rotatebox[origin=C]{180}{$\iota$}}
\newcommand{\iamp}{\rotatebox[origin=C]{180}{$\&$}}
\newcommand{\true}{\texttt{true}}
\newcommand{\false}{\texttt{false}}

\usepackage{url}

\title{Discernment is all you need%
\footnote{Dedicated to Marcelo Coniglio in his $60^{th}$ anniversary. Marcelo's enthusiasm for algebraic methods in logic sparked my own lasting fascination with equations, symmetries, and identity.}
}
\author{David Fuenmayor\\
    AI Systems Engineering \\
        University of Bamberg}
\date{}

\begin{document}

\maketitle

\begin{abstract}
We explore the expressive power of HOL, a system of higher-order logic, and its relationship to the simply-typed lambda calculus and Church's simple theory of types, arguing for the potential of HOL as a unifying logical framework, capable of encoding a broad range of logical systems, including modal and non-classical logics. Along the way, we emphasize the essential role of \textit{discernment}, the ability to tell things apart, as a language primitive; highlighting how it endows HOL with practical expressivity superpowers while elegantly enriching its theoretical properties.
\end{abstract}

\section{Funktionenkalkül and Logical Incompatibility}\label{intro}

On December 7, 1920, at a meeting of the G\"ottingen Mathematical Society, Moses Sch\"onfinkel (then working in Hilbert’s group) presented a talk titled \emph{Elemente der Logik}. Several years later, in March 1924, the contents of that talk were written up for publication in the \emph{Mathematische Annalen} as a paper entitled \emph{\"Uber die Bausteine der mathematischen Logik}.\footnote{The paper had to be finished by one of his colleagues (Heinrich Behmann), since Sch\"onfinkel had meanwhile left G\"ottingen---apparently for Moscow. In fact, few people, if any, knew about Sch\"onfinkel’s whereabouts after those years. Stephen Wolfram has recently carried out research on what may have happened; see \cite{wolfram2021did}.} That paper is a masterwork of clarity and depth and deserves to be summarized here in order to set the context.

Intending to lay out a minimalistic theory of atomic building blocks for mathematical logic in its most general conception, Sch\"onfinkel's paper \cite{schoenfinkel1924} starts (\S 1) by recalling the conceptual economy that Sheffer's stroke (interpreted as incompatibility, i.e.\ non-conjunction) brought to the propositional calculus \cite{Sheffer1913}, and proposes to extend this idea to the predicate calculus in the form of a variable-binding infix connective ``$\_\,|^{x}\,\_$'' representing incompatibility (\emph{Unvertr\"aglichkeit}). Thus
$$A(x)~|^x~B(x)$$
becomes interpreted as ``As and Bs are incompatible'', in the sense of A and B being mutually exclusive (as classes\footnote{For Schönfinkel, classes are special sorts of functions, namely \textit{propositional functions}, whose values are truth values.}), or, maybe more colloquially, as stating that if a thing is an A then it is certainly not a B (and vice versa). This can be seen as the strongest way of telling As and Bs apart: no-thing is ever both.

Sch\"onfinkel then digresses into the development of a variable-free \emph{Funktionenkalk\"ul}. This is the part that is often credited with introducing the concept of (what later became known as) ``combinators'',\footnote{It was around 1928 that Haskell Curry, then a graduate student at Princeton, (re)discovered Sch\"onfinkel's work and decided to move to G\"ottingen to complete his PhD dissertation (titled ``Grundlagen der kombinatorischen Logik'') within Hilbert's group. Initially, Curry had hoped to work with Sch\"onfinkel, but was informed by Pavel Alexandroff (then visiting Princeton) that Sch\"onfinkel had already left---apparently never to return (see \cite{wolfram2021did}). Nevertheless, Curry moved to G\"ottingen to work with Paul Bernays, one of the few scholars there who had worked with Sch\"onfinkel.} and which can be seen, with the benefit of hindsight, as an early and remarkably spare formalism for universal computation, well before Turing's seminal paper.

Sch\"onfinkel \textit{Funktionenkalkül} starts (\S 2) with an exposition of the underlying \textit{extensional} understanding of the concept of (mathematical) function as ``a correspondence between the elements of some domain of quantities, the argument domain, and those of a domain of function values", while emphasizing their higher-order nature ``permitting functions themselves to appear as argument values and also as function values''. It is here that Sch\"onfinkel presents the technique of arity reduction by iterated application, allowing for the systematic reduction of $n$-ary functions to unary ones, by treating an $n$-place function as returning, upon application to its first argument, an $(n{-}1)$-place function, so that expressions such as $Fxyz$ are parsed as $((Fx)y)z$.\footnote{This move has been dubbed ``Currying'' in modern literature, due to Curry's prominent use of it. Quine, who wrote the preface to the English translation of Sch\"onfinkel's paper \cite{schoenfinkel1924english}, states in it that this device was anticipated by Frege in \cite[\S 36]{Frege1893}.}

This technique naturally enables the introduction of some ``\textit{particular functions} of a very general nature'' (\S 3) by means of defining equations (used as rewriting rules), such as the identity function $I$ with $Ix = x$, the constancy function $C$ with $Cxy = x$, the interchange (``Ver\textit{t}auschung'') function $T$ with $Tfxy = fyx$, the composition (``\textit{Z}usammensetzung'') function $Z$ with $Zfgx = f(gx)$, and the fusion (``Ver\textit{s}chmelzung'') function $S$ with $Sfgx = fx(gx)$.\footnote{We have used Schönfinkel's original notation. These ``particular functions'' (excepting $S$) had been independently discovered by Curry in the late 1920's, who named them: $I$ (identity), $K$ (constancy), $C$ (interchange), $B$ (composition). Curry kept his own notation in his PhD thesis, and ever since, they have been known under those names, collectively, as ``combinators''.} Schönfinkel shows (\S 4) how they are all definable in terms of $S$ and $C$ only.

The previous conceptual reduction is applied (\S 5) to a ``special case'', namely, that of a ``calculus of logic in which the basic elements are individuals and the functions are propositional functions.'' This becomes possible by introducing ``an additional particular function, which is peculiar to this calculus'', and presents
$$U\,f\,g~=~f\,x~|^x~g\,x$$
as the defining equation of the incompatibility (``Unverträglichkeit'') function \textit{U}; commenting: \\
\textit{``It is a remarkable fact, now, that every formula of logic can be expressed by means
of our particular functions \emph{I}, \emph{C}, \emph{T}, \emph{Z}, \emph{S}, and \emph{U} alone, hence, in particular, by means solely of \emph{C}, \emph{S}, and \emph{U}.''}

Schönfinkel has essentially shown how to reduce logic to a calculus of functions (\emph{Funktionenkalkül}) extended with a primitive notion of incompatibility, and, even more impressively, to do so in what is known today as a point-free (i.e.\ variable-free, composition-driven) style; at the time epitomized by Peirce's 1870s algebra of absolute and relative terms. To get an idea of the scope of Schönfinkel’s contribution, we shall quote Quine in his preface of the English translation \cite{schoenfinkel1924english}:

\textit{``It was by letting functions admit
functions generally as arguments that
Schönfinkel was able to transcend the
bounds of the algebra of classes and relations and so to account completely for
quantifiers and their variables, as could
not be done within that algebra. The
same expedient carried him, we see, far
beyond the bounds of quantification
theory in turn; all set theory was his
province.''}

However, Sch\"onfinkel's presentation is still programmatic and leaves important pieces either implicit or only sketched: for instance, the formal status of the defining equations (as axioms vs.\ rewrite rules), the metatheory of the resulting equational calculus (confluence, normalization, consistency, etc.), and the systematic relationship between the variable-free function calculus and familiar logical formalisms are not worked out in anything like modern detail. These gaps were filled only later through the work of Curry (combinatory logic) and Church ($\lambda$-calculus, untyped and later typed), together with the contributions of their students and collaborators; see, e.g., \cite{CurryFeys1958,seldin2008church,hindley2008lambda} for detailed discussion.

Sch\"onfinkel's work provides a paradigmatic example of the idea of a calculus of functions (\textit{Funktionenkalkül}) as fundamental `plumbing' or `wiring' on top of which logical theories can be stated (and reasoned with) by adding the corresponding logical constants as primitive symbols (and axiomatizing them accordingly). Arguably, this idea had also been anticipated by Frege in his \textit{Begriffsschrift}, where logical inference is governed by an explicit function--argument structure and a precise notation for the representation of inferential steps. However, Frege's system remains tied to a specific logical vocabulary and proof calculus, and does not yet isolate the purely functional `plumbing' as an autonomous layer that can be reused and extended in a fully compositional way.

 Sch\"onfinkel's further insight, inspired by Sheffer's, was that a primitive notion of incompatibility (on top of a \textit{Funktionenkalk\"ul}) can play a foundational role for (classical) predicate logic.\footnote{Such an insight has been in the air for a long time (the author dimly recalls related remarks by medieval logicians). More recently, it has also been discussed by philosophers of logic, cf.~\cite{Brandom2008Incompatibility}, \cite{Peregrin2008LogicIncompatibility}.}
Indeed, the recent work of Coniglio and Toledo on \emph{logics of formal incompatibility} \cite{ConiglioToledo2023} illustrates this perspective in a contemporary setting, by taking incompatibility as a primitive \emph{binary} connective (written $\uparrow$) that localizes ``explosion''. Roughly: $A \uparrow B$, together with $A$ and $B$, trivializes deduction (i.e.\ anything follows---\textit{sequitur quodlibet}). Moreover, they show that the familiar ``consistency operator'' of the \textit{logics of formal inconsistency} (LFIs \cite{BookCC16}) arises as a special case of incompatibility, by setting $\circ A$ to be (or to behave like) $A \uparrow \neg A$, thereby exhibiting incompatibility as a genuine generalization of the notion of (in)consistency, and providing corresponding semantics and decision procedures for the resulting logical systems.

From a more abstract point of view, the notions of incompatibility at work in both Sch\"onfinkel's proposal (classical) and in the contemporary \textit{logics of formal incompatibility} of Coniglio and Toledo (paraconsistent) can be seen as a special case of a more general notion of \emph{discernment}. To avoid getting trapped in philosophical or etymological terrain, we shall now clarify that by \textit{discernment} we essentially mean the capacity to tell (different) things apart: the red from the green, the good from the bad, etc.\footnote{In some contexts, this ability can also be referred to as ``perception''. Without delving too deeply into psychological terrain, one could argue that it is, in a sense, an intuitive ``System 1''-type ability, which makes it reasonable to treat as primitive (i.e.\ undefinable) within a logical system. As is often noted in contemporary AI, this is precisely the kind of thing neural networks are very good at.} Incompatibility, on this reading, is simply one very harsh way of drawing such lines; it is not merely the failure of joint satisfiability, but a logical tool for marking distinctions: to say that $A$ and $B$ are incompatible is to treat them as mutually excluding options, and thereby to impose a discriminating structure on the space of assertions. In this sense, incompatibility can be taken as a primitive way of carving logical space, prior to (and more general than) negation or consistency. Sch\"onfinkel’s idea was that these basic discriminations could serve as foundational building blocks for logic (on top of a neutral \textit{Funktionenkalk\"ul}); Coniglio and Toledo's work, in turn, illustrate how this perspective can be made precise within modern proof-theoretic and semantic frameworks, in a paraconsistent setting.

Seen this way, incompatibility is best understood not as an isolated connective, but as one concrete manifestation of \emph{discernment}: the capacity of a logical system to articulate meaningful distinctions among its contents, on top of which notions of inference, consistency, and negation can then be built.
The next sections aim at introducing HOL, a classical higher-order logic, as a calculus of functions (the simply typed $\lambda$-calculus) extended with a primitive symbol of discernment (equality, or, if preferred, disequality). Before turning to that development, however, we still need to cover some more background.

\section{Logic in Typed $\lambda$-Calculi}\label{prelim}

What is now called the $\lambda$-calculus was introduced by Alonzo Church in the early 1930s as part of his program in the foundations of logic and mathematics. In fact, his initial formalism was not presented as a ``calculus'' but rather as ``a set of postulates for the foundation of logic'' \cite{Church1932} where ``rules of procedure'' (later called $\lambda$-conversion) were introduced together with logical constants, interrelated via ``formal postulates'', with the aim of serving as a foundation for logic and, by extension, mathematics.

It was this \emph{logical} system (i.e.\ $\lambda$-conversion plus the additional logical apparatus) that Church's students, S.\ Kleene and J.B.\ Rosser, famously showed to be inconsistent in 1935. Indeed, once one strips away the logical constants and logical postulates, the remaining formal `plumbing', the pure theory of $\lambda$-terms with conversion (aka.\ untyped $\lambda$-calculus), constitutes a consistent equational/reduction system in its own right that quickly proved to have extensive independent applications (notably as a framework for computation).

Moreover, it is worth emphasizing that Church's original presentation of $\lambda$-conversion was \textit{intensional},\footnote{Maybe the right wording would be \textit{hyperintensional}, depending who you ask---whether linguists or (modal) logicians.} since it did not build functional extensionality into the conversion rules; e.g., it did not employ anything like rule $\zeta$ ($M\,x = N\,x ~/~ M = N$) nor $\eta$-conversion ($\lambda x.\,M\,x = M$).\footnote{The usual proviso applies that $x$ does not appear free in $M$ nor $N$. In the early days, the $\eta$ and $\zeta$ conversion rules did not appear explicitly; they were isolated later on and shown equivalent to functional extensionality, cf.\ \cite[\S 3D]{CurryFeys1958} \cite[\S 7A]{hindley2008lambda}.} This stands in a contrast to Sch\"onfinkel’s \emph{Funktionenkalk\"ul}, where the extensional viewpoint is built into the intended interpretation and the way the defining equations are used (e.g.~his implicit use of rule $\zeta$).

It is fair to say, as is often stated in the literature, that the (untyped) $\lambda$-calculus is a \textit{calculus of functions}, conceived in their most idealized form: total, deterministic, and recursive. With the benefit of hindsight, it is perhaps unsurprising that early efforts by logicians (e.g., Church, Curry, and their students) to build \emph{systems of logic} on top of this conception of unrestricted, free-spirited functions ran into inconsistencies.

Here, so the story goes, the attempt to obtain consistent systems of logic based on such calculi of functions bifurcated into two main directions: the ``simple theory of types'' of Church, and the (primarily untyped) ``combinatory logic'' pioneered by Curry. In this paper we shall focus on the former. We refer the reader to \cite{seldin2008church} for a discussion of the historical development of both combinatory logic and $\lambda$-calculus.

\subsection{The Simple Theory of Types}\label{sec-stt}

In the present writing, a starring role is played by Church's ``simple theory of types'' (STT), so we shall briefly recall its presentation below (see \cite{sep-type-theory-church} for a proper discussion).

The set $\T$ of STT's types is inductively generated, starting from a set of type constants $\T\C = \{\bool,\iota\}$ by applying the binary type constructor $\ar$ (written as a right-associative infix operator; see Fig.~\ref{fig:HOL-grammar})
\begin{figure}
	\normalsize	\centering
\begin{align*}
\alpha,\, \beta \;&::=\;\; \tau \in \T\C  \;\;\; | \;\;\; \alpha \,\ar\, \beta \\
s,\, t \,\;&::=\;\; c_\alpha \in \Const \;\; | \;\; x_\alpha \in \Var \;\; | \;\; \left(\lambda x_\alpha.\, s_\beta\right)_{\alpha\ar\beta}
\;\; | \;\; \left(s_{\alpha\ar\beta} \; t_\alpha\right)_\beta
\end{align*}
\caption{Informal grammar for STT types and terms.}
\label{fig:HOL-grammar}
\end{figure}

For instance, $\bool$,  $\bool\,{\ar}\,\bool$, $\iota\,{\ar}\,\iota$ and $\iota\,{\ar}\,\iota\,{\ar}\,\bool$ are all well-formed types. STT can be seamlessly extended by assuming a larger set of type constants $\T\C \supseteq \{\bool,\iota\}$ as needed for applications. It is worth noting that the axioms of STT (see \cite{sep-type-theory-church}) place enough constraints such that, semantically speaking, the type $\bool$ ends up having exactly two different inhabitants (modulo interderivability), and thus they are called ``Boolean'' truth-values. By contrast, the type $\iota$ is left unconstrained (often intuitively interpreted as the type of ``individuals''). Note that Church's original formulation of STT includes functional extensionality as an explicit but optional axiom, which we adopt here.

STT \emph{terms} are inductively defined starting from a collection of typed constant symbols ($\Const$) and typed variable symbols ($\Var$) using the constructors \emph{function abstraction} and \emph{function application}, by obeying the corresponding type constraints (see Fig.~\ref{fig:HOL-grammar}).
Type subscripts and parentheses are usually omitted to improve readability, if obvious from the context or irrelevant. In contrast to predicate logic, in STT there is no distinction between \textit{terms} and \textit{formulas}. We observe that STT terms of type `$\bool$' are customarily referred to as ``propositions'', and sometimes as ``formulas'' too. 

In the original STT presentation by Church \cite{church1940formulation}, and also the one by Henkin \cite{henkin1950completeness}, an infinite number of logical constants are introduced (to be constrained by axiom schemata). We gather them in a set $\Const $.\footnote{Church (and Henkin) actually employed the symbols $N$, $A$, and $\Pi$ for \texttt{not}, \texttt{or}, and $\forall$, respectively. We cannot follow Church in all of his notational choices in this paper.}
$$\Const = \{\texttt{not}_{\bool\ar\bool}, \texttt{or}_{\bool\ar\bool\ar\bool}\} \cup \{\forall^{\alpha}_{(\alpha\ar\bool)\ar\bool} ~|~ \alpha\in\T\} \cup \{ \iiota^{\alpha}_{(\alpha\ar\bool)\ar\alpha} ~|~ \alpha\in\T\}$$
The families $\forall^\alpha$ and $\iiota^\alpha$ represent quantifier and definite description predicates for all types $\alpha\in\T$.\footnote{In fact, both Church and Henkin seamlessly employ schematic type meta-variables ($\alpha$, $\beta$, etc.) when listing their constants and axiom schemata (via a less subtle abuse of notation than ours), somehow anticipating the need for schematic type-polymorphism, as featured in the STLC (to be discussed in \S\ref{sec-stlc}).}
The other Boolean connectives are introduced as abbreviations in the expected fashion. The customary infix notation for binary connectives ($\vee$, $\rightarrow$, etc.) is also introduced.

It is worth reminding the reader that the term constructor \emph{function abstraction} (aka.\ $\lambda$-abstraction) is the only variable-binding construct available in the STT. The customary binder notation for quantifiers and definite descriptions is introduced as an abbreviation (`syntactic sugar'):
\begin{align*}
\forall^\alpha x_{\alpha}.\,s_\bool &\;:\equiv\; \forall^\alpha_{(\alpha \ar \bool)\ar \bool}\,(\lambda x_{\alpha}.\,s_{\bool}) &\text{for each}~\alpha \in \T \\
\iiota^\alpha x_{\alpha}.\,s_\bool &\;:\equiv\; \iiota^\alpha_{(\alpha \ar \bool)\ar \alpha}\,(\lambda x_{\alpha}.\,s_{\bool}) &\text{for each}~\alpha \in \T
\end{align*}

Church introduces equality in STT as a family of abbreviations, following Leibniz's principle of \textit{identity of indiscernibles}:\footnote{Observe that this formulation can be shown equivalent to the one employing double-implication (by contraposition, noting that quantifiers range over \textit{all} predicates, including their negations). As it happens, this definition, often called ``Leibniz-equality'', while intended to denote the identity relation, might fail to do so in some ``non-standard'' models. We will elaborate on this issue later on (see \S\ref{sec-q}), since first we need to introduce some basics of STT semantics.}
$$ \Q^\alpha_{\alpha\ar\alpha\ar\bool} :\equiv \lambda x_\alpha. \lambda y_\alpha. \forall^{\alpha\ar\bool} f_{\alpha\ar\bool}. (f~x \rightarrow f~y) ~~~\text{for each}~\alpha \in \T$$

Using the previously defined equality symbol $\Q$ (with infix notation $=$), Church proposes two axiom schemata for the constants $\iiota^\alpha$, namely, the ``axioms of descriptions''
$$f_{\alpha\ar\bool}\,x_\alpha \rightarrow (\forall y_\alpha .~f~y \rightarrow x = y) \rightarrow f\,(\iiota_{(\alpha\ar\bool)\ar\alpha}^\alpha\,f)$$
and the ``axioms of choice''
$$f_{\alpha\ar\bool}\,x_\alpha \rightarrow f\,(\iiota_{(\alpha\ar\bool)\ar\alpha}^\alpha\,f)$$
noting that the former clearly follow from the latter. Church's intention was to allow the use of $\iiota^\alpha$ as \textit{definite description} operators without necessarily committing to ``choice'' principles (in which case $\iiota^\alpha$ effectively becomes a choice/selection operator à la Hilbert's $\varepsilon$). Indeed, the later schema directly entails the following (more stereotypical) formulation of the ``axiom of choice'' in mathematics, whereby every indexed family of non-empty sets (i.e.\ left-total relation) $R_{\beta\ar\gamma\ar\bool}$ has a choice function $g_{\beta\ar\gamma}$ (corresponding to $\lambda x_\beta.\,\iiota_{(\gamma\ar\bool)\ar\gamma}^\gamma (R\,x))$:
$$\forall R_{\beta\ar\gamma\ar\bool}.\,(\forall x_\beta.\,\exists y_\gamma.\,R\,x\,y) \rightarrow (\exists g_{\beta\ar\gamma}.\,\forall x_\beta.\,R\,x\,(g\,x))$$

Moreover, Church mentions that the terms $\iiota_{(\alpha\ar\bool)\ar\alpha}^\alpha$ (intended as \textit{definite descriptions}\footnote{Unfortunately, Church does not elaborate on this point. Later on Henkin \cite{henkin1963theory} and Andrews \cite[\S 53--$5309^\gamma$]{andrews2002} verify this claim when $\iiota^\alpha$ are axiomatized via the ``axiom of descriptions'' above. However, when $\iiota^\alpha$ are interpreted/axiomatized as choice operators (i.e.\ Hilbert's $\varepsilon$), the possibility of such a reduction collides with later theoretical results (see \cite{Yasuhara1975} \cite[\S1.3.5]{sep-type-theory-church}).}) are definable inductively for any functional type $\alpha$ of form $\beta\ar\gamma$. Thus, if $\iiota_{(\gamma\ar\bool)\ar\gamma}^\gamma$ is already defined, we can define:	
$$\iiota^{\beta{\ar}\gamma}_{((\beta{\ar}\gamma){\ar}\bool){\ar}\beta{\ar}\gamma} :\equiv~ \lambda h_{(\beta{\ar}\gamma){\ar}\bool}.\,\lambda x_\beta.\, \iiota^\gamma y_\gamma. \, \exists f_{\beta{\ar}\gamma}. \,h\;f \land y =^\gamma f\,x$$
    
So it seems like we only need to take definite descriptions for the base type constants ($\{\iota,\bool\}$ for STT) as primitives. As it happens, only $\iiota^{\iota}_{(\iota{\ar}\bool){\ar}\iota}$ is required, since we can in fact define:\footnote{This was noted later on by Henkin \cite{henkin1963theory}. Andrews \cite{andrews1972a} also provides several such definitions.}
$$\iiota^{\bool}_{(\bool{\ar}\bool){\ar}\bool} :\equiv \Q^{\bool{\ar}\bool}~(\lambda x_\bool.\,x)$$

It is worth noting that, in keeping with the Zeitgeist, Church did not introduce a ``semantics'' for STT. This occurred ten years later at the hands of Leon Henkin (another of Church's students), to whom we owe the first completeness proof for a calculus of STT.
For our current purposes,\footnote{For a properly detailed presentation of STT semantics, we refer the reader to \cite{sep-type-theory-church,andrews2002}.} it shall suffice to mention that in Henkin's \textit{general models} \cite{henkin1950completeness} (a paradigmatic set-based approach to STT semantics) each model is based upon a \textit{frame}: a collection $\{\DD_\alpha\}_{\alpha \in \T}$ of non-empty sets, called \textit{domains} (for each type $\alpha$). $\DD_\bool$ is chosen as a two-element set, say $\{T,F\}$, in order to make STT classical, whereas $\DD_\iota$ may have arbitrarily many elements (aka.\ ``individuals''). As expected, the set $\DD_{\alpha \ar \beta}$ consists of functions with domain $\DD_\alpha$ and codomain $\DD_\beta$. In so-called ``standard'' models \cite{sep-type-theory-church}, these domain sets are assumed to be \textit{full}, i.e.~they contain \textit{all} functions from $\DD_\alpha$ to $\DD_\beta$.

Famously, Henkin \cite{henkin1950completeness} enlarged the class of STT models to also include ``non-standard'' ones, where domains $\DD_{\alpha \ar \beta}$ are not full, yet contain all functions `named' by an STT term, by means of a suitably defined \textit{denotation function} ($|\cdot|$) which interprets each term $s_\alpha$ as an element $|s_\alpha|$ of $\DD_\alpha$ (its \textit{denotation}). As expected, a \textit{denotation function} worthy of its name must respect the intended semantics of STT as a \textit{logic of functions}. Thus, $| s_{\alpha\ar\beta} \; t_\alpha | \in \DD_\beta$ denotes the value of the function $|s_{\alpha\ar\beta}| \in \DD_{\alpha\ar\beta}$ when applied to $|t_\alpha| \in \DD_\alpha$, and $|\lambda x_\alpha.\, s_\beta| \in \DD_{\alpha\ar\beta}$ denotes the corresponding function from $\DD_\alpha$ to  $\DD_\beta$ (see e.g.~\cite[\S2]{sep-type-theory-church} for details on what this means in terms of substitutions and the like). Denotations of terms are defined inductively, in the expected way, starting with the term constants, which in STT are \texttt{not}, \texttt{or}, and $\forall$.

\subsection{The Simply-Typed Lambda Calculus} \label{sec-stlc}

It should be noted that what computer scientists nowadays call the ``simply-typed lambda calculus'' (STLC) has arisen as an \textit{a posteriori} distillation of Church's STT. In its basic version, it can be defined similarly to STT but with all type- and term-constants removed (as well as their constraining axioms), leaving only the functional plumbing/wiring. Seen from this perspective, it might also be fair to say that STLC is a calculus (but not yet a logic) of \textit{typed} functions, i.e., functions with a clearly delimited domain and codomain. Moreover, we shall adopt here an \textit{extensional} understanding of functions (e.g.~by assuming $\eta$-conversion). After all, we want the STLC to play the role of a \textit{Funktionenkalk\"ul}.

Lacking type constants (such as $\iota$ and $\bool$), modern presentations of the STLC employ schematic type variables (cf. \textit{rank-1} or `prenex' type polymorphism in the literature). This is very useful in practice, for instance, it makes best sense to define a function: 
$$\texttt{applyTwice} :\equiv \lambda f.~\lambda x.~f~(f~x)$$
that applies a given function twice, as having the type $$(\alpha~\ar~\alpha)~\ar~(\alpha~\ar~\alpha)$$
parametric in the type $\alpha$ (think of it as a schema variable). Essentially, this means that such a schematic-type-polymorphic function only needs to be defined once, and it will work uniformly for any type $\alpha$ chosen by the caller.\footnote{This feature of type systems turned out to be of such tremendous importance for programming that even mainstream object-oriented languages like Java and Go ended up implementing it (as so-called ``generics'') after years of programmers' protests.}

In practical applications (e.g., in programming languages), the STLC plays the role of a formal framework that can be extended (or instantiated) by adding domain-specific type- and term-constants (as well as constraining axioms) as needed. Seen from this perspective, the STT \textit{could} be understood, quite anachronistically of course, as an instantiation of the STLC, as applied to the domain of logic.\footnote{This modern understanding of the STLC is difficult to trace back to the work of some particular individual(s), so Church is usually credited due to his intellectual influence.}

There is a great deal of literature on the expressive strengths (and weaknesses) of the STLC, so we shall not delve into it here. Instead, we restrict ourselves to revisiting a traditional classroom exercise for $\lambda$-calculus students; namely, defining terms that behave like the Boolean constants \texttt{true} and \texttt{false}, and then using these to define terms that behave like pairs, conditionals, and Boolean-like connectives such as conjunction and disjunction (which we pedantically denote with an asterisk as a reminder that they are not quite the real thing). We start with
\[
\begin{aligned}
\texttt{true}^* &\;:\equiv\; \lambda x.\,\lambda y.\,x\\
\texttt{false}^* &\;:\equiv\; \lambda x.\,\lambda y.\,y
\end{aligned}
\]

Next, we encode pairs as
\[
\langle x, y \rangle \;\;:\equiv\;  \lambda f.\,f\,x\,y
\]
so we can extract the first and second components using projections
\[
\begin{aligned}
\pi_1 &\;:\equiv\; \lambda p.\,p~\texttt{true}^*\\
\pi_2 &\;:\equiv\; \lambda p.\,p~\texttt{false}^*
\end{aligned}
\]
thus getting\footnote{The notation $\cong$ means that both expressions reduce to the same normal form according to the $\beta$- and $\eta$-conversion rules of the $\lambda$-calculus.}
\[
\pi_1~\langle x, y \rangle \cong x ~~~~\text{resp.}~~~
\pi_2~\langle x, y \rangle \cong y.
\]

As it happens, a conditional operator can be defined in terms of pairs
\[
\texttt{if}^* \;:\equiv\; \lambda c.\,\lambda t.\,\lambda e.\,\langle t,e\rangle~c \;\;\;\cong\; \lambda c.\,\lambda t.\,\lambda e.\,c~t~e
\]
so that
\[
\begin{aligned}
\texttt{if}^*~\texttt{true}^*~t~e 
&\;\;\;\cong\;\;\; \langle t,e\rangle~\texttt{true}^* &\cong\;\; \pi_1~\langle t,e\rangle &\;\;\cong\; t\\
\texttt{if}^*~\texttt{false}^*~t~e
&\;\;\;\cong\;\;\; \langle t,e\rangle~\texttt{false}^* &\cong\;\; \pi_2~\langle t,e\rangle &\;\;\cong\; e
\end{aligned}
\]

Using the previous definitions we can define some `Boolean-like' connectives

\[
\begin{aligned}
\texttt{not}^* &\;:\equiv\; \lambda b.\,\texttt{if}^*~b~\texttt{false}^*~\texttt{true}^*  &\cong\;\;&\lambda b.\,b \;\texttt{false}^* \;\texttt{true}^*\\
\texttt{and}^* &\;:\equiv\; \lambda b_1.\,\lambda b_2.\,\texttt{if}^*~b_1~b_2~\texttt{false}^* &\cong\;\;&\lambda b_1.\,\lambda b_2.~b_1~b_2~\texttt{false}^* \\
\texttt{or}^*  &\;:\equiv\; \lambda b_1.\,\lambda b_2.\,\texttt{if}^*~b_1~\texttt{true}^*~b_2 &\cong\;\;&\lambda b_1.\,\lambda b_2.~b_1~\texttt{true}^*~b_2
\end{aligned}
\]

so that the expected (trivial) equivalences are satisfied, for instance:
\[
\begin{aligned}
\texttt{not}^*~\texttt{true}^* &\cong \texttt{false}^*
&\qquad
\texttt{not}^*~\texttt{false}^* &\cong \texttt{true}^*\\
\texttt{and}^*~\texttt{true}^*~b &\cong b
&\qquad
\texttt{and}^*~\texttt{false}^*~b &\cong \texttt{false}^*\\
\texttt{or}^*~\texttt{true}^*~b &\cong \texttt{true}^*
&\qquad
\texttt{or}^*~\texttt{false}^*~b &\cong b
\end{aligned}
\]

Hence the definitions above seem to capture the usual operational behavior of pairs, conditionals and basic logical operations \textit{within the \(\lambda\)-calculus}, and can be adapted straightforwardly to the simply typed setting (STLC) by adding the appropriate type annotations.

However, as useful as the above encodings are (e.g.\ in the theory of programming languages), to the practicing logician they have the feeling of being not more than an \textit{ersatz}. Things quickly start to break down the moment we try more serious logic stuff. For instance, consider the following claim, corresponding to the De Morgan law:
\[
\texttt{not}^*\,(\texttt{and}^*~a~b) \cong \texttt{or}^*~(\texttt{not}^*~a)~(\texttt{not}^*~b).
\]

The above does not hold unrestrictedly under arbitrary values for the variables $a$ and $b$. There is the implicit constraint that they must be instantiated only with values ($\lambda$-expressions) from the set $\{\texttt{true}^*, \texttt{false}^*\}$ for this to hold. On a related note (to the dismay of the unaware reader), we shall note that
\[
\langle \pi_1~p\,,\pi_2~p \rangle \;\ncong\; p .
\]

Further red flags arise when we attempt to assign types to expressions featuring the above encodings. While the given definitions are simple enough to have well-formed types in the STLC, strange things happen with more complex expressions. For example, when instantiating the previous formulation of De Morgan's law with the pseudo-Boolean terms $\texttt{true}^*$ and $\texttt{false}^*$, they end up having \textit{necessarily distinct} types (and quite convoluted ones, indeed). Similarly, we would intuitively expect the $\lambda$-expressions $\texttt{and}^*$ and $\texttt{or}^*$ above to have the same type (after all, they are both `dual' binary Boolean connectives). Well, they don't.

\subsection{Equality as Primitive and the System $\Q_0$} \label{sec-q}

So, back to the notion of a \textit{logic} of functions (like STT), in contrast to a mere \textit{calculus} of functions (like STLC), we see that the idea of equality (identity) as the primitive fundamental logical connective (notion) is certainly not new. It was apparently in the air at the beginning of the last century, and being actively worked on by Ramsey, possibly inspired by the philosophical views of Wittgenstein, as we can infer from this passage from Ramsey's ``The Foundations of Mathematics''\footnote{Quoted in \cite{andrews2014bit}, where Andrews recalls Henkin's development of a formulation of type theory based on equality, and the significance of this contribution.}

\textit{``The preceding and other considerations led Wittgenstein to the view that mathematics does not consist of tautologies, but of what he called `equations', for which I would prefer to substitute `identities' \dots (It) is interesting to see whether a theory of mathematics could not be constructed with identities for its foundation. I have spent a lot of time developing such a theory, and found it
was faced with what seemed to me insuperable difficulties.''}

 As noted by Leon Henkin in his paper ``Identity as a logical primitive'' \cite{henkin1975identity}, Alfred Tarski, as early as 1923 \cite{tarski1923}, had already noted that, in the context of higher-order logic (which he calls ``la Logistique'' heavily influenced by the \textit{Principia}), one can define propositional connectives in terms of logical equivalence and quantifiers. We take the opportunity to extract from that paper by Henkin our working definition of identity:\footnote{The terms ``identity'' and ``equality'' are used in the literature almost interchangeably. We prefer to think of ``equality'' as a connective (syntax) and ``identity'' as the denoted relational concept (semantics).}

\textit{``By the relation of identity we mean that binary relation which holds between any object and itself, and which fails to hold between any two distinct objects.''} \cite[p.\,31]{henkin1975identity}

In that paper (to which we refer the reader interested in a deeper treatment), Henkin credits Quine \cite{quine1937logic, quine1956unification} with having shown how both quantifiers and connectives can be defined in terms of equality and the abstraction operator $\lambda$ in the context of Church's type theory. Henkin later came up (independently \cite{henkin1963theory}) with an analogous formulation and provided the corresponding axiom system, which also benefited from further simplifications by Peter Andrews (another student of Church) \cite{andrews1963reduction,andrews2014bit}. The main definitions are presented as follows (with other connectives introduced as abbreviations in the expected way):
\[
\begin{aligned}
\;\true_\bool \;\;\;\;\;\;\;&:\equiv~ (\lambda x_\bool.\,x) =^{\bool{\ar}\bool} (\lambda x_\bool.\,x) \\
\;\false_\bool \;\;\;\;\;\;\;&:\equiv~ (\lambda x_\bool.\,x) =^{\bool{\ar}\bool} (\lambda x_\bool.\,\true) \\
\;\texttt{not}_{\bool{\ar}\bool} \;\;\;\;\;\;&:\equiv~ (\lambda x_\bool.\,\false =^\bool x) \\
\forall^\alpha_{(\alpha{\ar}\bool){\ar}\bool} &:\equiv~ (\lambda P_{\alpha{\ar}\bool}.\,P =^{\alpha{\ar}\bool} (\lambda x_\alpha.\,\true)) \\
\texttt{and}_{\bool{\ar}\bool{\ar}\bool} \;\;\,&:\equiv~ (\lambda x_\bool.\,\lambda y_\bool.\,(\lambda f_{\bool\ar\bool}. f\,x =^\bool y) =^{(\bool\ar\bool)\ar\bool} (\lambda g_{\bool\ar\bool}.\,g~\true))
\end{aligned}
\]

Let us try to understand the above encoding. The constant \texttt{true} can be an arbitrary tautology, like claiming that something (e.g.~the identity function) is self-identical. The constant \texttt{false} can be an arbitrary absurdity or impossibility, like claiming of two different things that they are identical. Negation becomes self-explanatory once we observe that the axioms constrain the semantical domain for type `$\bool$' to contain exactly two inhabitants (the denotations of \texttt{true} and \texttt{false}). The universal quantifier basically tells of a given predicate whether it holds always true for any input. Conjunction is admittedly not very intuitive; maybe this prompted Andrews to introduce the following alternative definition (see \cite[\S1.1]{Andrews1965TTT} and also \cite[\S1.4]{sep-type-theory-church}, \cite[\S51]{andrews2002}):
$$\texttt{and}_{\bool{\ar}\bool{\ar}\bool} \;\;\,:\equiv~ (\lambda x_\bool.\,\lambda y_\bool.\,(\lambda f_{\bool{\ar}\bool{\ar}\bool}.\,f\,x\,y) =^{(\bool{\ar}\bool{\ar}\bool){\ar}\bool} (\lambda f_{\bool{\ar}\bool{\ar}\bool}.\,f\,\true\,\true)) $$

The above is basically an encoding of the corresponding truth table (recalling the ordered pair encoding from STLC in \S\ref{sec-stlc}, it can be written as: $\lambda x. \lambda y.~\langle x,y\rangle = \langle \true,\true\rangle$). 

This simplified system featuring equality as the sole (primitive) logical constant would later give rise to the system $\Q_0$ in Andrews' doctoral dissertation \cite{Andrews1965TTT}, where it is introduced as the foundational subsystem for his transfinite type system $\Q$. Ever since, $\Q_0$ has become a standard foundational reference for higher-order automated reasoning and interactive theorem proving.\footnote{Notably as the basis of Andrews' TPS/ETPS line of provers and closely related to the HOL family; see \cite{andrews2002} for a detailed textbook presentation and metatheoretic development.}

Finally, we note that $\Q_0$ (following STT) features a family of functions $\iiota_{(\alpha\ar\bool)\ar\alpha}^\alpha$ axiomatized to behave as inverses of equality, i.e., for any $y_\alpha$ we have\footnote{This axiom schema is equivalent to the previously discussed ``axioms of descriptions'' in Church's STT. Andrews also presents ``axioms of choice'' as optional extensions to his system.}
$$(\iiota^\alpha x_\alpha.~ y_\alpha=^\alpha x) ~\cong~ \iiota_{(\alpha\ar\bool)\ar\alpha}^\alpha (\lambda x_\alpha.~y_\alpha =^\alpha x) ~\cong~ \iiota^\alpha (\Q^\alpha y_\alpha) ~=^\alpha~ y_\alpha$$
Recalling from \S\ref{sec-stt}, we observe that, from the above family, only $\iiota^\iota_{(\iota\ar\bool)\ar\iota}$ needs to be taken as primitive. In fact, being the inverse operation to $\Q^\alpha$, $\iota^\alpha$ can be seen as yet another instance of discernment.\footnote{By extension, the same can be said of Hilbert's $\varepsilon$. Choice is a facet of discernment.}

It is worth mentioning that Andrews' postulation of equality (resp.~identity) as a logical primitive goes beyond any minimalistic or aesthetic appeal, as it ultimately becomes a mathematical byproduct (resp.~necessity) of Henkin’s notion of general models (see discussion in \S\ref{sec-completeness}). More specifically, Andrews constructs in \cite{andrews1972general} a general (and non-standard) model in which functional extensionality fails to hold.\footnote{His is in fact a quite concrete model with $| \DD_{\iota}| = 3$. Adopting equality as primitive has the added benefit of eliminating finite non-standard models \cite[\S 54]{andrews2002}, and thus any generated finite (counter)model is necessarily ``standard''.} His diagnostic is basically that sets (i.e.~elements of $\DD_{\alpha\ar\bool}$) in such a (non-standard) model may be so sparse that the denotation of Leibniz-equality is not the actual identity relation. In fact, Andrew's result entails that Henkin's soundness theorem for STT \cite[Theorem 2]{henkin1950completeness} is ``technically incorrect'' (in Andrews' words), yet Henkin's completeness result remains unaffected. As Andrews aptly points out \cite{andrews1972general}, the issue is more conceptual than mathematical, and it is resolved by modifying the definition of general models to include the requirement that $\DD_{\alpha \ar \bool}$ contains all singleton sets.

Moreover, Andrews further argues in \cite{andrews1972general} that, since the identity relation has now made its way into the semantic definition of STT, the most natural next step is to introduce it too into the syntax as a primitive connective, namely equality, while possibly defining the rest in terms of it, as illustrated in $\Q_0$.\footnote{Of course, an alternative approach would be to sacrifice extensionality instead; cf.~\cite{J6}, who systematically explore generalizations of general models, among them models where functional extensionality and so-called Boolean extensionality ($|\DD_\bool| = 2$) fail.} In fact, later work in the context of automated theorem proving in HOL does indeed introduce equality as primitive---for good practical reasons (see \cite{B5} and references therein).

\section{HOL as an Universal (Meta-)Language}

\subsection{What is HOL?}
It is often said, half-jokingly, that formal logic sits somewhat uncomfortably between the philosophy and computer science chairs---with only the occasional excursion into mathematics. Perhaps the same could be said of the acronym HOL, which, at this point, we can only describe as standing for ``higher-order logic''. But what exactly is HOL?

In philosophy (and sometimes mathematics departments), people tend to think of HOL, usually under the influence of set theory, quite literally as ``$\alpha$-order logic'', for some ordinal $\alpha$, meaning that it features quantifiers ranging over $\alpha{-}1$-order predicates, all the way down to second-order logic (which features quantifiers ranging over good-old first-order predicates). This view has the benefit of reusing an arguably well known conceptual framework (predicates, quantifiers, relations, etc.).

In computer science, HOL is employed as an \textit{umbrella term} for a family of \textit{classical} higher-order logics, built upon STLC, and featuring some kind of weak type polymorphism, like the schematic (`prenex') one discussed previously. More specifically, in theorem proving \cite{B5}, HOL is employed to refer to logical languages used in several mathematical proof assistants (Isabelle/HOL, HOL-Light, HOL4, etc.) descending from the venerable \textit{HOL} system \cite{gordon1985hol}. In this text we will abstract away the particular implementation-specific differences between HOL `flavors'. Thus, we present HOL as an \textit{idealized} logical (family of) system(s) which feature the best of STLC and $\Q_0$ worlds: They support schematic type-polymorphism and use it to introduce a primitive term constant $\Q$ having type $\alpha \ar \alpha \ar \bool$ parametric on an arbitrary type $\alpha$.\footnote{Of course, pragmatic implementation details diverge among systems. For instance Gordon's HOL \cite{gordon1985hol} introduces implication as an additional primitive connective (instead of defining it using equality). Others also feature primitive quantifiers. These design choices are mostly made to boost performance (e.g.\ calculi rules often feature implication and quantifiers directly).}

Thus, HOL is not only a \textit{higher-order} functional calculus (i.e.\ it allows for functions that take other functions as arguments and/or return functions) but also becomes a \textit{higher-order logic}: we can define quantifiers that range over predicate variables (of type $\alpha \ar \bool$) and, more generally, arbitrary function variables (type $\alpha\ar \beta$).

When it comes to expressivity, higher-order logic (HOL) has been attacked from both sides: too much (e.g.~being undecidable) and too little (having a too `simple' type system in which even type inference is decidable). We shall take this as good evidence that HOL indeed reaches an expressivity sweet spot. To be fair, FOL proponents can also concoct an analogous sweet-spot argument, and maybe they are right too. After all, HOL (with an appropriate semantics\footnote{One with respect to which sound and complete calculi exist. Such a semantics being e.g.~Henkin's general models \cite{henkin1950completeness}, under which HOL satisfies first-order model-theoretical properties (compactness, L\"owenheim-Skolem, etc.). We refer the reader to \cite{Andreka2014} and \cite{manzano2017identity,manzano2018identity} for an extended discussion.}) and many-sorted FOL (and thus vanilla FOL) are equally expressive from a Lindstr\"om's theorem perspective \cite{manzano2017identity, farmer2008seven}. Yet, the pragmatics of formalizing and proving in FOL and HOL couldn't feel more different. In terms of practical applications, HOL seems to have the upper hand, as witnessed by the fact that all full-fledged mathematical proof assistants are based on (some variation of) HOL\footnote{Systems like ACL2, Mizar and Metamath being the kind of exceptions that confirm the rule, as they add additional syntactic HOL-like sugar: ACL2 provides functional syntax and type-like guards, while Mizar also emulates `types' and higher-order functions via schemes (among others). Moreover, Mizar and Metamath are based on set-theory (``HOL in wolf's clothing''). Still these systems are not as popular as their HOL-based cousins.}.

In a sense, the main thesis of the present writing essentially claims that HOL (understood as STLC plus equality/disequality) is ``all you need''. Still, as with love, people sometimes crave additional expressivity. For them,\footnote{Many of whom are researchers working on the formalization of mathematics using dependently-typed proof assistants like PVS, Coq, Lean, etc.} a natural extension of HOL featuring dependent types has recently been developed. This system is fittingly called \textit{dependently-typed HOL} (DHOL).
Nonetheless, HOL is expressive enough to embed DHOL \cite{C94}. This possibility is not surprising if we recall the well-known analogous encoding of HOL into (many-sorted) FOL.

\subsection{On Completeness}\label{sec-completeness}

It is not uncommon to hear in logic folklore that HOL is ``incomplete''. Of course, without further elaboration, such a statement sounds as puzzling as it is meaningless. People who say this (when they are not quoting hearsay) are usually referring to the notion of ``Gödel incompleteness'' (aptly dubbed ``essential incompleteness'' by Andrews \cite{andrews2002}), or perhaps more appropriately, ``incompletability''. Indeed, this metaphysical property is proudly worn by HOL as a badge of honor: it testifies that HOL is capable of expressing the most interesting kinds of problems.\footnote{Philosophical speculations around ``G\"odel incompleteness'' abound, and it is not our aim to give them further platform. We refer the interested reader to \cite{own20GoedelParaconsistency} for own speculations, or \cite{franzen2005} for a sober perspective.}

Another, more down-to-earth notion of (in)completeness is a relational property of a given calculus with respect to an intended semantics (e.g.\ a class of mathematical structures). When showing that a particular calculus is incomplete (with respect to a semantics), we provide evidence of a blind spot—usually a `bug' in the proof procedure.\footnote{In automated reasoning, it is not uncommon to consciously add optimizations that break completeness of calculi in exchange for improved efficiency. The lost completeness will ideally make it to the backlog as `technical debt'.} This is a problem that can, in principle, be fixed—e.g., by adding more axioms or rules of inference. 
By contrast, ``Gödel incompleteness'' is, by its very nature, unfixable. It is therefore not a `bug', but rather a feature.

Sound and complete STT- and HOL-calculi have existed since at least Henkin's seminal 1950 paper \cite{henkin1950completeness}. They are plentiful, so we refer the reader to \cite{B5} for a somewhat dated but very insightful survey. Henkin's main insight consists in discovering an appropriate semantics for interpreting STT (and thus HOL), so that complete calculi exist (cf.~\cite{Andreka2014} for a discussion). The reader shall note that this semantics is often qualified with the adjective ``general'' in the literature (e.g. Henkin's \textit{general models}), whereas the traditional (inadequate) one often gets the qualification ``standard''.\footnote{As in ``standard movie'', ``standard excuse'', or ``standard deontic logic''}

Recall from our previous discussion (\S\ref{sec-stt}) that, broadly speaking, in HOL's ``standard'' semantics, the domain sets $\DD_{\alpha \ar \beta}$ (which interpret functional types) are \textit{full}, i.e.\ they are supposed to contain \textit{all} functions with domain $\DD_\alpha$ and codomain $\DD_{\beta}$; and thus they are completely determined by the choice of domains for the base type constants (e.g.~$\DD_\iota$ for STT). This convenient property comes at the cost of metaphysical opacity (e.g.~what exactly is in $\DD_{\iota\ar\bool}$ for infinite $\DD_\iota$?).

By contrast, in general models, the domains $\DD_{\alpha \ar \beta}$ are not necessarily \textit{full}, but still contain enough elements (functions) to guarantee that any term $s_{\alpha \ar \beta}$ has a denotation. In other words, the only requirement for the domains $\DD_{\alpha \ar \beta}$ is that they contain, at least, all \textit{nameable} functions.\footnote{This can be equated  to a \textit{nominalist} position, as opposed to an arguably \textit{platonist} one in the case of ``standard'' models; cf.~\cite{henkin1953some, fernandes2022henkin}.} On the flip side, the family of domains $\DD_{\tau}$ (called a ``frame'') is no longer determined by the choice of domains for base types. So far, it seems this is the conceptual bullet we have to bite in order to define an appropriate HOL semantics, that is, one that allows for the possibility of a complete calculus.

But maybe we can have our cake and eat it too. Towards the end of the aforementioned article \cite{henkin1975identity}, Henkin discusses ``Boolean models'' for HOL,\footnote{In a footnote, Henkin mentions that this work was originally presented at a symposium on \textit{Theory of Models} held in Berkeley in 1963, but the corresponding paper was never published (no reason given). He does not provide any proofs, which he assures will be given ``in a forthcoming paper''—which the present author has not been able to locate.} in which, he claims, this convenient situation arises. In $\B$-\textit{models} (as Henkin names them) the domain $\DD_\bool$ can be an arbitrary \textit{complete Boolean algebra} $\B$ (not necessarily the $2$-valued one from before). More specifically, he introduces $\B$-models for a relational version of STT, as being completely determined by choosing a non-empty domain $\DD_\iota$ and a function $\mathbf{q}$ which associates to each pair of individuals $x$, $y$ in $\DD_\iota$ a value $\mathbf{q}(x,y)$ in the complete Boolean algebra $\B$. As Henkin further explains, ``this function $\mathbf{q}$ may be thought of as measuring `the degree of equality' of $x$ and $y$''.\footnote{As expected, this function $\mathbf{q}$ must satisfy special equivalence-like properties: It must be commutative, and it should validate $\mathbf{q}(x,x) = 1$, as well as $\mathbf{q}(x,y) \cdot \mathbf{q}(x,z) \leq \mathbf{q}(y,z)$ for all $x,y,z \in \DD_\iota$ (corresponding to variants of symmetry, reflexivity and Euclideanness).}

Once $\DD_\iota$ and $\mathbf{q}$ are chosen, all other domains are determined: We take as $\DD_\bool$ the carrier set of $\B$, and interpret $\Q_{(\bool,\bool)}$ (equality among propositions)\footnote{In this relational variant of STT, a relational type $(\alpha_1, \dots, \alpha_n)$ corresponds to an $n$-ary relation among $\DD_{\alpha_i}$ for $i \leq n$. It can be paraphrased (via Currying) as the functional type $\alpha_1 \ar \dots \ar \alpha_n \ar \bool$ in STT.} as the operation $\Leftrightarrow$ of double-implication in $\B$; $\Q_{(\iota,\iota)}$ denotes $\mathbf{q}$; finally, $\Q_{(\beta,\beta)}$ for $\beta = (\alpha_1, \dots \alpha_n)$ denotes the function $\mathbf{q}^\beta(r,s) = \bigwedge(r(x_1, \dots, x_n) \Leftrightarrow s(x_1, \dots, x_n))$ for all $x_i \in \DD_{\alpha_i}$. As in standard models of STT, $\DD_{(\alpha_1, \dots, \alpha_n)}$ is defined to be the \textit{full} domain of relations among domain sets $\DD_{\alpha_i}$ for $i \leq n$. Furthermore, Henkin teases in  \cite{henkin1975identity} with the existence of a stronger completeness result: ``There is a \textit{specific} complete Boolean algebra $\B^*$, such that for every formula $A_\bool$ we have $\models_{\B^*} A_\bool$ if, and only if, $\vdash A_\bool$.''

We finish this discussion by noting that standard models are subsumed under general models, and thus HOL formulas proven valid with respect to the general semantics are also valid in the standard sense.  Moreover, we shall observe that (classes of only) non-standard models cannot be characterized via HOL formulas \cite[\S55]{andrews2002}. As a consequence, it is not clear how (if at all) the results delivered by HOL calculi might ever differ from the ``standard'' ones. In the words of Andrews \cite[p.\,255]{andrews2002}: \textit{``One who speaks the language of [HOL] cannot tell whether he lives in a standard or nonstandard world, even if he can answer all the questions he can ask.''}

\subsection{Modal and Non-classical Logics in HOL}

A technique called \textit{shallow semantical embeddings} \cite{J41,J23} has been developed to encode (quantified) modal and non-classical logics into HOL as a meta-language, in such a way that object-logical formulas correspond to HOL terms. This is realized by directly encoding in HOL the truth conditions (semantics) for object-logical connectives as syntactic abbreviations (definitions), essentially in the same way as they appear in textbooks. 

In a sense, this is not much different from what logicians have been doing since the invention of model-theoretical semantics, where traditional logic expositions can be seen as `embedding' the semantics of object logics into natural language (e.g., English plus mathematical shorthand) as a meta-language. The difference now is that the meta-language, HOL, is itself a formal system of logic.\footnote{This is, in fact, the idealized situation Tarski envisioned in his seminal paper \cite{Tarski1933Truth}. He had to fall back on natural language for technical reasons (he was clearly ahead of his time!).}

This way we can carry out semi-automated  object-logical reasoning by translation into a formal meta-logic (HOL) for which good automation support exists \cite{B5}. As argued, e.g., in \cite{J41}, such a \textit{shallow} embedding allows us to reuse state-of-the-art automated theorem provers and model generators for reasoning with (and about) many different sorts of non-classical logical systems in a very efficient way (avoiding, e.g., inductive definitions and proofs as in \textit{deep} embeddings).\footnote{We refer the reader to \cite{DeepShallow} for a discussion of the differences (and similarities) between deep and shallow embeddings.}

The idea of employing HOL as meta-logic to encode \textit{quantified} non-classical logics has been exemplarily discussed by Benzm\"uller \& Paulson in \cite{J23} for the case of normal multi-modal logics featuring first-order and propositional quantifiers. Their approach draws upon the ``propositions as sets of worlds'' paradigm from modal logic, by adding the twist of encoding sets as their characteristic functions, i.e., as total functions with a (2-valued) Boolean codomain, in such a way that a set $S$ becomes encoded as the function $s$ of type $\alpha\ar\bool$ such that for any $x$ of type $\alpha$ we have that $x \in S$ iff $(s~x) = \true$. As has been discussed in  \cite{own23SemanticalInvestigations}, this basically corresponds to Stone-type representations of \textit{Boolean Algebras with Operators} \cite{BAO} as algebras of sets. For the sake of illustration, we recall the corresponding definitions in HOL (see \cite{own23SemanticalInvestigations} for details), employing $\iota$ for the type of `worlds' (and thus $\iota\ar\bool$ for propositions).

\[
\begin{aligned}
\U_{\iota\ar\bool} \;\;&:\equiv~ \lambda w_\iota.~\true \\
\emptyset_{\iota\ar\bool} \;\;&:\equiv~ \lambda w_\iota.~\false \\
\--_{(\iota\ar\bool)\ar{\iota\ar\bool}} \;&:\equiv~ \lambda P_{\iota\ar\bool}.~\lambda w_\iota.~\neg (P~w) \\
\cup_{(\iota\ar\bool)\ar(\iota\ar\bool)\ar{\iota\ar\bool}}\,&:\equiv~ \lambda P_{\iota\ar\bool}.~\lambda Q_{\iota\ar\bool}.~\lambda w_\iota.~ (P~w) \vee (Q~w)\\
\cap_{(\iota\ar\bool)\ar(\iota\ar\bool)\ar{\iota\ar\bool}}\,&:\equiv~ \lambda P_{\iota\ar\bool}.~\lambda Q_{\iota\ar\bool}.~\lambda w_\iota.~ (P~w) \wedge (Q~w)\\
\Rightarrow_{(\iota\ar\bool)\ar(\iota\ar\bool)\ar{\iota\ar\bool}}\,&:\equiv~ \lambda P_{\iota\ar\bool}.~ \lambda Q_{\iota\ar\bool}.~ \lambda w_\iota.~ (P~w) \rightarrow (Q~w) \\
\leftharpoondown_{(\iota\ar\bool)\ar(\iota\ar\bool)\ar{\iota\ar\bool}}\,&:\equiv~ \lambda P_{\iota\ar\bool}.~ \lambda Q_{\iota\ar\bool}.~ \lambda w_\iota.~ (P~w) \wedge \neg (Q~w) \\
\bigcup\nolimits_{((\iota\ar\bool)\ar\bool)\ar{\iota\ar\bool}} \;&:\equiv~ \lambda S_{(\iota\ar\bool)\ar\bool}.~\lambda w_\iota.\,~\exists P_{\iota\ar\bool}.~ (S~P) \wedge (P~w) \\
\bigcap\nolimits_{((\iota\ar\bool)\ar\bool)\ar{\iota\ar\bool}} \;&:\equiv~ \lambda S_{(\iota\ar\bool)\ar\bool}.~\lambda w_\iota.\,~\forall P_{\iota\ar\bool}.~ (S~P) \rightarrow (P~w) \\
\subseteq_{(\iota\ar\bool)\ar(\iota\ar\bool)\ar\bool}\,&:\equiv~ \lambda P_{\iota\ar\bool}.~ \lambda Q_{\iota\ar\bool}.~ \forall w_\iota.~ (P~w) \rightarrow (Q~w) \\
\end{aligned}
\]

In this way, after interpreting Boolean connectives ($\dot\top, \dot\bot, \dot\vee, \dot\wedge, \dot\rightarrow$, etc.)\footnote{We shall sometimes add a dot to differentiate between object-logical and HOL connectives.} as their corresponding set operations, modalities become easily encoded as shown below---not unlike the well-known \textit{standard translation} of modal logic into first-order logic. \footnote{Note that they are further parameterized with an argument $R$ of type $\iota \ar \iota \ar \bool$ in the role of accessibility relation. Hence the present encoding corresponds in fact to multi-modal logic.}
\[
\begin{aligned}
\square_{(\iota\ar\iota\ar\bool)\ar(\iota\ar\bool)\ar{\iota\ar\bool}} \;&:\equiv~ \lambda R_{\iota\ar\iota\ar\bool}.~\lambda P_{\iota\ar\bool}.~\lambda w_\iota.\,~\forall v_\iota.~ (R~w~v) \rightarrow (P~v) \\
\lozenge_{(\iota\ar\iota\ar\bool)\ar(\iota\ar\bool)\ar{\iota\ar\bool}} \;&:\equiv~ \lambda R_{\iota\ar\iota\ar\bool}.~\lambda P_{\iota\ar\bool}.~\lambda w_\iota.\,~\exists v_\iota.~ (R~w~v) \wedge (P~v) 
\end{aligned}
\]

On top of this, Benzm\"uller \& Paulson \cite{J23} have shown how to encode object-logical quantifiers (first-order and propositional) by lifting the meta-logical ones.
\[
\begin{aligned}
\dot\forall_{(\alpha\ar\iota\ar\bool)\ar\iota\ar\bool} \;&:\equiv~ \lambda \varphi_{\alpha\ar\iota\ar\bool}.~ \lambda w_\iota.\, \forall x_\alpha.~ (\varphi~x~w)\\
\dot\exists_{(\alpha\ar\iota\ar\bool)\ar\iota\ar\bool} \;&:\equiv~ \lambda \varphi_{\alpha\ar\iota\ar\bool}.~ \lambda w_\iota.\, \exists x_\alpha.~ (\varphi~x~w)
\end{aligned}
\]

Binder notation can also be introduced in the usual way, where $\dot\forall x.\,\varphi$ stands for $\dot\forall (\lambda x.\,\varphi)$ and $\dot\exists x.\,\varphi$ stands for $\dot\exists (\lambda x.\,\varphi)$.

Following Benzm\"uller \& Paulson work, a series of papers have been brought forward to substantiate the claim that HOL can truly serve as a universal meta-logic for many different systems of non-classical logic. For instance, we have paracomplete \& paraconsistent logic \cite{own23SemanticalInvestigations}, deontic logics \cite{own19HarnessingHOLComplexEthicalTheories,J68}, epistemic and dynamic logic \cite{J58}, conditional and defeasible logic \cite{C37,J68}, among many others specimens in the non-classical logics zoo.

For the sake of illustration, we provide below a selection (necessarily out of context) of some non-classical connectives, just to give the reader a sense of how their HOL encodings look. We refer to the corresponding papers for context and discussion. Also note that we omit pedantic type-annotations---for readability, and because, in practice, they can be automatically inferred (e.g.\ by Isabelle/HOL):

Constant- resp.\ variable-domain universal quantifiers, parameterized with a set $E$ resp.\ relation $\E$ representing `existing' objects (see \cite{own17AutomatingEmendations}, and \cite[\S 3.8]{own23SemanticalInvestigations}).
\[
\begin{aligned}
\dot\forall^{\texttt{con}} \;&:\equiv~ \lambda E.~\lambda \varphi.~ \lambda w.\, \forall x.~ (E~x) \rightarrow (\varphi~x~w)\\
\dot\forall^{\texttt{var}} \;&:\equiv~ \lambda \E.~\lambda \varphi.~ \lambda w.\, \forall x.~ (\E~x~w) \rightarrow (\varphi~x~w)
\end{aligned}
\]

Boolean-based LFIs' \cite{own22LFIsCarnielliConiglio} paraconsistent negation and its corresponding consistency recovery operator, assuming a unary set-operation $\B_{(\alpha\ar\bool)\ar\alpha\ar\bool}$ as (axiomatized) constant representing the topological notion of ``border'' (see \cite{own23SemanticalInvestigations}).
\[
\begin{aligned}
\dot\neg \;&:\equiv~ \lambda P.\,\lambda x.~{\neg}(P~x) \vee (\B~P~x) &\text{i.e.}~~\lambda P.~{-}P \cup (\B~P)\\
\circ \;&:\equiv~ \lambda P.\,\lambda x.~{\neg}(\B~P~x) &\text{i.e.}~~ \lambda P.~{-}(\B~P)
\end{aligned}
\]

Fusion and its residual implications (e.g.~as in Lambek calculus and also some relevance logics) via ``Routley-Meyer-style'' semantics,\footnote{We refer the reader to \cite{beall2012ternary} for explanations about this kind of semantics in the context of relevance logic, and analogously to \cite[\S 2]{Moortgat2014} for the Lambek calculus. Isabelle/HOL encodings of these systems can be consulted in \cite[\texttt{examples/substructural\_logics}]{LogicBricks}.} assuming a ternary relation $R_{\alpha\ar\alpha\ar\alpha\ar\bool}$ as constant (duly axiomatized).

\[
\begin{aligned}
\otimes \;&:\equiv~ \lambda P.\,\lambda Q.\,\lambda z.\,\exists x.\,\exists y.\, (P~x) \wedge (Q~y) \wedge (R~x~y~z)\\
\text{\textbackslash} \;&:\equiv~ \lambda P.\,\lambda Q.\,\lambda y.\,\forall x.\,\forall z.\, (R~x~y~z) \rightarrow (P~x) \rightarrow (Q~z) \\
/ \;&:\equiv~ \lambda P.\,\lambda Q.\,\lambda x.\,\forall y.\,\forall z.\, (R~x~y~z) \rightarrow (Q~y) \rightarrow (P~z)
\end{aligned}
\]

In the above examples, propositions denoted sets (more specifically, elements of a complete Boolean algebra of sets\footnote{Even more specifically, in HOL (and $\Q_0$) such a Boolean algebra is also atomic, with singletons (i.e.~denotations of $\Q\,x$ for any $x$) playing the role of atoms.}), and operations denoted set operations (e.g. `left-image' as in the unary $\lozenge$ or the binary $\otimes$), but in principle, propositions can denote elements in any algebraic structure as long as it admits a suitable ordering for defining ``local'' (aka.\ ``degree-preserving'' \cite{Font2007}) consequence\footnote{Otherwise, we can still put our creativity to work and try to come up with an adequate ``global consequence'' definition, e.g., using generalized matrices and the like; cf.~\cite{Wojcicki1970, Font2007} for a discussion. The literature on many-valued logics has quite some good examples.} and features suitably interrelated operations (e.g.~adjunction/residuation conditions enabling things like modus ponens, deduction theorem, and the like). 

Related to the modal ``propositions as sets (of states)'' paradigm, we can also conceive of an analogous ``propositions as relations (on states)'' one. Understanding propositions as sets of transitions between states offers an enriching, dynamic perspective. In order to do this, we shall first encode an algebra of relations \textit{qua} set-valued functions. Semantically speaking, a relation with source $\DD_\alpha$ and target $\DD_\beta$ is represented as the function $R \in \DD_{\alpha \ar \beta \ar \bool}$ (via Currying). From this perspective, relations inherit the structures of both sets and functions, and enrich them in manifold ways.

Literature on relational algebra informally classifies relational connectives into two groups: a so-called ``Boolean'' (aka.\ additive) and a ``Peircean'' (aka.\ multiplicative) structure. When seen as set-valued functions, we can see that the former is in fact inherited from sets, while the latter comes from their generalized (i.e., partial and non-deterministic) functional structure. We start with the former, by doing essentially the same `intensionalization' move as for sets before (and again omitting type annotations):
\[
\begin{aligned}
\U^r \;\;&:\equiv~ \lambda w.~\lambda v.~\true \\
\emptyset^r \;\;&:\equiv~ \lambda w.~\lambda v.~\false \\
\--^r \;&:\equiv~ \lambda R.~\lambda w.~\lambda v.~\neg (R~w~v) ~~&(\text{also notation:}~R^{-})\\
\cup^r\,&:\equiv~ \lambda R_1.~\lambda R_2.~\lambda w.~\lambda v.~ (R_1~w~v) \vee (R_2~w~v)\\
\cap^r\,&:\equiv~ \lambda R_1.~\lambda R_2.~\lambda w.~\lambda v.~ (R_1~w~v) \wedge (R_2~w~v)\\
\Rightarrow^r\,&:\equiv~ \lambda R_1.~ \lambda R_2.~ \lambda w.~\lambda v.~ (R_1~w~v) \rightarrow (R_2~w~v) \\
\leftharpoondown^r\,&:\equiv~ \lambda R_1.~ \lambda R_2.~ \lambda w.~\lambda v.~ (R_1~w~v) \wedge \neg (R_2~w~v) \\
{\bigcup}^r \;&:\equiv~ \lambda S.~\lambda w.~\lambda v.\,~\exists R.~ (S~R) \wedge (R~w~v) \\
{\bigcap}^r \;&:\equiv~ \lambda S.~\lambda w.~\lambda v.\,~\forall R.~ (S~R) \rightarrow (R~w~v) \\
\subseteq^r\,&:\equiv~ \lambda R_1.~ \lambda R_2.~ \forall w.~\forall v.~ (R_1~w~v) \rightarrow (R_2~w~v) 
\end{aligned}
\]
The algebraic structure above is essentially the same as that for sets (complete atomic Boolean algebra). Thus, in the context of HOL, relations can seamlessly be employed in the same roles sets are (unsurprisingly, since they are, after all, isomorphic to sets of pairs). More interesting is their ``Peircean'' structure:\footnote{We refer to \cite{LogicBricks} for a discussion, and several applications of these relational notions.}
\[
\begin{aligned}
\smallsmile \;&:\equiv~ \lambda R.~\lambda w.~\lambda v.~(R~v~w) ~~&(\text{also notation:}~R^\smallsmile)\\
\smallfrown \;&:\equiv~ \lambda R.~\lambda w.~\lambda v.~\neg (R~v~w) ~~&(\text{also notation:}~R^\smallfrown)\\
\textbf{;} \,&:\equiv~ \lambda R_1.~\lambda R_2.~\lambda w.~\lambda v.~\exists u.~(R_1~w~u) \wedge (R_2~u~v)\\
\dagger \,&:\equiv~ \lambda R_1.~\lambda R_2.~\lambda w.~\lambda v.~\forall u.~(R_1~w~u) \vee (R_2~u~v)\\
\vartriangleright \,&:\equiv~ \lambda R_1.~\lambda R_2.~\lambda w.~\lambda v.~\forall u.~(R_1~u~w) \rightarrow (R_2~u~v)\\
\blacktriangleleft \,&:\equiv~ \lambda R_1.~\lambda R_2.~\lambda w.~\lambda v.~\forall u.~(R_2~v~u) \rightarrow (R_1~w~u)
\end{aligned}
\]

From the operations on relations above, the first two are unary: $\smallsmile$ is transposition (aka.\ ``converse'', ``reverse'', etc.), $\smallfrown$ is cotransposition (converse-of-complement/complement-of-converse); and the rest are binary operations (written as infix notation): \textbf{;} is composition, $\dagger$ is dual-composition, and $\vartriangleright$ resp.~$\blacktriangleleft$ are residuals `on the right' resp.~`on the left' (wrt.~composition).\footnote{Note that: $R \dagger S = (R^{-} ~\textbf{;}~ S^{-})^{-}$. We speak of residuals, in the sense that: $R~\textbf{;}~S \subseteq^r T$ iff $S \subseteq^r R \vartriangleright T$ iff $R \subseteq^r T \blacktriangleleft S$. Moreover, we have that $R \vartriangleright S = R^\smallfrown \dagger S$ and $R \blacktriangleleft S = R \dagger S^\smallfrown$.}

As we can see, HOL is not only a \textit{logic of functions}, but also hides a very powerful \textit{logic of relations}. As an example, we can use the previous insights to provide a shallow semantical embedding for \textit{cyclic  linear logic},\footnote{This linear logic variant was introduced in~\cite{Desharnais1997}. See \cite[\texttt{examples/substructural\_logics}]{LogicBricks} for the corresponding Isabelle/HOL encoding.} basically by associating object-logical connectives to relational operations as follows:\footnote{Note that $\D$ corresponds to disequality, which can be encoded as: $\lambda x.\lambda y. \neg (x = y)$.}
\[
\begin{aligned}
\textbf{1} \;&:\equiv~ \Q &\text{~~~~~~and~~~~~~} \textbf{0} \;&:\equiv~ \emptyset^r\\
\bot \;&:\equiv~ \D &\text{~~~~~~and~~~~~~} \top \;&:\equiv~ \U^r\\
\otimes \;&:\equiv~ \textbf{;} &\text{~~~~~~and~~~~~~} \oplus \;&:\equiv~ \cup^r\\
\iamp \;&:\equiv~ \dagger &\text{~~~~~~and~~~~~~} \& \;&:\equiv~ \cap^r\\
\multimap \;&:\equiv~ \vartriangleright &\text{~~~~~~~and~~~~~} {\bullet}\hspace{-2pt}{-} \;&:\equiv~ \blacktriangleleft\\
\end{aligned}
\]

Finally, ``exponentials'' are also encoded as unary operations on relations:
$$\textbf{!} \;:\equiv~ \lambda R.~R ~\cap^r~ \Q \text{~~~~~~and~~~~~~} \textbf{?} \;:\equiv~ \lambda R.~R ~\cup^r~ \D$$

\subsection{A Recipe for Encoding Logics in HOL}
We provide below an informal recipe we have been following to encode logics in HOL. We have worked mostly with the Isabelle/HOL proof assistant, but the approach transfers seamlessly to other systems (HOL-Light, Coq, Lean, etc.).

\begin{itemize}
	\item Make a decision: What shall propositions denote—sets, relations, categories, finite field elements or polynomials? This choice constrains many future decisions.
	\item Get acquainted with the corresponding background mathematical theory for the semantics, which is typically part of some library for the proof assistant in question (and usually includes algebra, set theory, topology, category theory, etc.).
	\item State the signature for object-logical connectives and map arities to types. For instance, if your propositions have type $\alpha$, your candidate concepts for unary and binary connectives might have types of the form $\alpha \ar \alpha$ and $\alpha \ar \alpha \ar \alpha$, respectively—or maybe $\beta \ar \alpha \ar \alpha$ and $\beta \ar \alpha \ar \alpha \ar \alpha$ in case they take an additional argument of type $\beta$ (e.g., as in relation-indexed modalities in multi-modal logics).
	\begin{itemize}
		\item Use these types to search for (or generate) candidate concepts in the background mathematical theory to interpret the connectives in question (e.g.~$\lozenge$ resp.~$\otimes$ are interpreted as unary resp.~binary left-image set-operators for a binary resp.~ternary relation).
		\item Decide on the type of entailment: value- vs.\ degree-preserving (aka.\ global vs.\ local) \cite{Font2007}. Provide notation (e.g., in modal logics local consequence ends up being syntactic sugar for subset ordering).
		\item Encode semantic axioms (e.g., frame conditions in modal logic).
		\item Carry out some tests. They can be `positive' like verifying axioms/rules and known theorems (using proof automation, e.g.~\textit{Sledgehammer} \cite{blanchette2016hammering}). They can also be `negative' like finding counterexamples to expected non-theorems (using a model generator, e.g., \textit{Nitpick} \cite{blanchette2010nitpick}).
		\item Go back to some previous step as required, and iterate \dots\footnote{As an instance of the ``problem of formalization'', the logic embedding process is interpretative, and thus iterative and virtuously circular: it shall converge to a state of ``reflective equilibrium''. See previous work on \textit{computational hermeneutics} \cite{own19CHApproachConceptualExplicitation, own18ACaseStudyCH} for a discussion.}
	\end{itemize}
\end{itemize}

After surveying the previous, more or less sophisticated, encodings of modal and other ``non-classical'' logics in HOL, we might wonder where HOL's enhanced expressivity comes from, as compared with the poor-man logic encoding hacks of the STLC (as discussed in~\S\ref{sec-stlc}). 
As we have seen, the secret sauce is \textit{discernment}: all you need.

\section{Discernment and Duality: Prospects}

We have previously mentioned that by \textit{discernment} we essentially mean the capacity to tell (different) things apart: the red from the green, the good from the bad, etc. It can be mathematically modelled, quite naturally, as a relation, in two dual ways, depending on where we want to place the focus. They are:
\begin{itemize}
    \item{\textbf{Equality/identity:}} The relation denoted by $\Q$ (infix $=$).
    \item \textbf{Disequality/difference:} The \textit{dual} relation denoted by $\D$ (infix $\neq$).
\end{itemize}

The concept of (classical) negation serves as a mediator in this duality: ${a \neq b}$ iff ${\neg(a = b)}$, and ${a = b}$ iff ${\neg(a \neq b)}$. Perhaps this reveals its most essential characteristic: a bridge between sameness and difference.

We saw before how $\Q$ added to STLC, thus giving rise to HOL, is enough to define all other logical connectives and quantifiers. We also briefly mentioned the \textit{shallow semantical embeddings} technique \cite{J23,J41}, whereby HOL can embed all kinds of logics as fragments.\footnote{The universalistic claim is, in fact, as strong as it sounds. We are not aware of any systems of non-classical logic that cannot be embedded in HOL—sometimes in a quite natural and elegant way. We are, of course, very much interested in obtaining concrete evidence of the limits of this approach.}
Hence, after unfolding definitions in the mentioned HOL encodings, all we are left with is $\Q$ plus the functional STLC wiring.\footnote{Some Isabelle/HOL encodings introduce additional constants that are constrained by axioms, e.g., as in the case of accessibility relations in modal logics (but their role is actually that of free variables). Still, $\Q$ (resp.~$\D$) remains the sole logical constant in the language.} Identity is all you need.

We also saw how adding $\Q$ (or its evil twin $\D$) as a term-constant to the STLC (to obtain HOL) presupposes the existence of  a ``Boolean'' type constant `$\bool$' too. That this type constant is interpreted as a domain set $\DD_\bool$ having exactly two elements, does not only follow tradition, but also corresponds to an intuitive ``all or nothing'' understanding of identity (difference). However, we saw from Henkin's efforts \cite{henkin1975identity} described above, that a more general, algebraic understanding of ``Boolean'' is also possible in this context, and thus identity (difference) can also be understood in ``shades of grey'', even rendering arguably more elegant completeness results.

Thus, encoding logical connectives using discernment as primitive, either as identity or as difference, gives rise to the pair of dual translations below:
\\
\\
\textbf{Via positiva:}
\[
\begin{aligned}
\;\true \;\;&:\equiv~ \Q = \Q \\
\;\false \;\;&:\equiv~ (\lambda x.\,\true) = (\lambda x.\,x) \\
\;\texttt{not} \;&:\equiv~ (\lambda x.\,x = \false) \\
\texttt{and} \;\,&:\equiv~ (\lambda x.\,\lambda y.\,(\lambda f.\,f\,x\,y) = (\lambda f.\,f\,\true\,\true)) \\
\forall^\alpha \;&:\equiv~ (\lambda P.\,P = (\lambda x_\alpha.\,\true)) \\
\dots
\end{aligned}
\]
\textbf{Via negativa:}
\[
\begin{aligned}
\;\false \;\;&:\equiv~ \D \neq \D \\
\;\true \;\;&:\equiv~ (\lambda x.\,\false) \neq (\lambda x.\,x)  \\
\;\texttt{not} \;&:\equiv~ (\lambda s.\,s \neq \true) \\
\texttt{or} \;\,&:\equiv~ (\lambda s.\,\lambda t.\,(\lambda f.\,f\,s\,t) \neq (\lambda f.\,f\,\false\,\false)) \\
\exists^\alpha \;&:\equiv~ (\lambda P.\,P \neq (\lambda x_\alpha.\,\false)) \\
\dots
\end{aligned}
\]

In fact, the notion of \textit{duality}, as it appears throughout logic and lattice theory (among others), can often be traced back to the fundamental duality (mediated by negation) between $\Q$ and $\D$. If you have an axiom, definition, or conjecture, and wish to find its `dual', you can reduce it to STLC plus discernment, and then simply switch $\Q$s with $\D$s.\footnote{Advertisement: This insight has been instrumental in constructing the \textit{Combinatory Logic Bricks} Isabelle/HOL library \cite{LogicBricks}, aimed at providing a one-stop shop for formalized mathematical notions commonly used in modal and non-classical logics.} It is interesting to consider whether duality can serve as a lighthouse, guiding future investigations into lattice-valued models for HOL or Boolean-valued models for (co-)set theory (e.g., with coinduction treated as a first-class citizen). A first interesting exercise would be to reconstruct the results (and proofs) claimed by Henkin for his Boolean(-valued) HOL models,\footnote{As an interesting side remark, as stated in~\cite{henkin1975identity}, Henkin developed these models in 1963---the same year in which Cohen published his famous paper introducing the ``forcing'' method.  Boolean-valued models were developed by Solovay and Scott a couple of years later, and reportedly, there was a manuscript by Scott circulating in 1966 titled ``Boolean-valued Models for Higher-Order Logic'' (possibly part of the infamously unpublished ``Scott–Solovay paper''). Moreover, Scott mentions in his foreword to \cite{bell1977boolean} that ``in September of 1951, in a paper of Alonzo Church \cite{church1953-mex} delivered at the Mexican Scientific Conference, a suggestion for Boolean-valued models of type theory had already been made.''} in terms of both $\Q$ and $\D$ (e.g.~with $\D_{\bool\ar\bool\ar\bool}$ interpreted as symmetric-difference/xor) and all definitions suitably dualized.

A further interesting line of work concerns HOL models based on Boolean algebras (or lattices) ``with operators'' \cite{BAO}. Previous work~\cite{own23SemanticalInvestigations} flirts with this idea in the context of \textit{Topological Boolean Algebras} (TBAs), featuring topologically motivated unary operators like closure, interior, frontier, border, and the like. 
More concretely, in the area of algebraic fuzzy logic, there has been some model-theoretical work worth highlighting, such as Novák's \textit{Fuzzy Type Theory} (FTT)~\cite{novak2005fuzzy}, in which the domain set $\DD_{\bool}$ corresponds to an \textit{IMTL-algebra}, namely, a residuated lattice (satisfying prelinearity and double negation) extended by the \textit{Baaz delta} operation~\cite{Hajek1998Metamathematics}. In that spirit, Běhounek~\cite{behounek2016} introduces ``a minimalistic many-valued theory of types'' ($TT_0$), which also generalizes STT by allowing arbitrary algebras of truth-values, providing a foundational, modular framework into which Nov\'ak's FTT fits as an (algebraically richer) instance focused on graded truth and fuzzy equality.

In a similar, though less algebraic, spirit, the work of Kohlhase and Scheja~\cite{kohlhase1999higher} builds upon Smullyan's~\cite{Smullyan1963} \textit{abstract consistency} technique (as later extended to STT by Andrews~\cite{andrews1971resolution} and to multi-valued first-order logics by Carnielli~\cite{carnielli1987systematization}) to provide a systematic treatment of completeness in resolution calculi for many-valued variants of STT (in which $\DD_{\bool}$ is a finite set). 
A very interesting possibility in this regard adds an algebraic touch to such finite-many-valued approaches by leveraging insights from the theory of finite (aka.\ Galois) fields; in particular, by treating an $n$-valued algebra of truth-values as a field of order $n$ (provided $n$ is a prime-power), or more simply as the field of integers modulo~$n$ (when $n$ is prime). Duality insights can help ensure a smooth transition from the classical case ($\mathbb{F}_2$) to higher ones and guide the choice of definitions throughout the process. Consequently, (meta-)logical investigations can profit from the arsenal of algebraic tools associated with (Galois) fields, like equational solving techniques, Gr\"obner bases, Lagrange interpolation, etc.,\footnote{For example, formulas (and also connectives) can be represented as polynomials, and (at least in the propositional case) satisfiability/provability amounts to checking whether (sets of) polynomials have (common) roots.} as well as from nice uniqueness results (like polynomials having unique normal forms, unique factorizations, etc.).
These insights have been applied to propositional (see e.g.~\cite{chazarain1991multi,carnielli2005polynomial,carnielli2014non, agudelo2016polynomial}), modal \cite{agudelo2011polynomial,Agudelo2019modal}, and first-order \cite{carnielli2018reconciling} logics, but not yet to higher-order ones.\footnote{More advertising: Some preliminary Isabelle/HOL  sources formalizing this approach (and interfacing with \textit{computer algebra systems}) are available in GitHub \cite{algebraic-theorem-proving}.}

Finally, there is a substantial body of work on extensions of STT with partial functions, which we shall not survey here, as it is quite extensive and not directly aligned with our fuzzy-algebraic perspective.\footnote{Interested readers may consult recent work such as~\cite{farmer2023simple} and~\cite{Manzano2023Hybrid}, and the references therein.} 

In a sense, it should come as no surprise to claim that (vanilla) HOL, as discussed here, can shallowly embed all the higher-order systems mentioned above as fragments as well, much in the same spirit as previously presented embeddings of intensional-modal HOL~\cite{own17AutomatingEmendations,ownAFPTypesTableaux} or dependently-typed HOL~\cite{C94}. 
Whether it can do so in a practically effective way (e.g., for the purposes of automated reasoning), however, may require less a mathematical proof and more a technological (and sociological) one.\footnote{Indeed, some empirical evidence for the effectiveness of semantic embeddings in HOL has been shown recently by Steen et al.\ \cite{SteenEtAl2025}, using Steen's logic embedding tool \cite{SteenLET2022}, in the context of a comparative study of automated reasoners for quantified modal logics. As it turns out, despite legitimate concerns about the indirection overhead of HOL embeddings, it was shown that the latter can often outperform native reasoning systems.}

\bibliographystyle{abbrv}
\bibliography{all}

\end{document}